\documentclass{amsart}
\usepackage{xcolor}
\usepackage{amsfonts,amssymb,amsmath,amsthm}
\newtheorem{theorem}{Theorem}[section]
\newtheorem{proposition}[theorem]{Proposition}
\newtheorem{corollary}[theorem]{Corollary}
\newtheorem{example}[theorem]{Example}
\newtheorem{lemma}[theorem]{Lemma}

\newtheorem{remark}[theorem]{Remark}

\DeclareMathOperator{\lin}{sp}

\DeclareMathOperator{\rank}{rank}
\begin{document}
\title{Representation of compact operators between Banach spaces}
\author{G. Ramesh}
\address{G. Ramesh, Department of Mathematics, IIT Hyderabad, Kandi, Sangareddy, Telangana- 502284, India.}
\email{rameshg@math.iith.ac.in}
\author{M. Veena Sangeetha}
\address{M. Veena Sangeetha, Department of Mathematics, IIT Hyderabad, Kandi, Sangareddy, Telangana- 502284, India.}
\email{veena176@gmail.com}
\author{Shanola S. Sequeira}
\address{Shanola S. Sequeira, Department of Mathematics, IIT Hyderabad, Kandi, Sangareddy, Telangana- 502284, India.}
\email{shanolasequeira@gmail.com}
\keywords{Banach space, compact operator, norm attaining operator, spectral theorem, duality mapping}
\subjclass{46B50, 47A67, 15A18}
\maketitle

\begin{abstract}
	In this article, we give a representation for compact operators acting between reflexive Banach spaces, which generalizes the representation given by Edmunds et al. for compact operators between reflexive Banach spaces with strictly convex duals. Further, we give a representation for operators on Banach spaces that are comparable to compact normal operators on Hilbert spaces and illustrate our result with an example.
\end{abstract}

\section{Introduction}
Let $X$ be a Banach space over the field $\mathbb{F}\in\{\mathbb{R}, \mathbb{C}\}$. We denote the unit sphere and the dual of $X$ by $S_X$ and $X^*$, respectively. If $X$ and $Y$ are two Banach spaces over the same field, we denote the Banach space of all bounded linear operators from $X$ into $Y$  by $\mathcal{B}(X,Y)$ and the subspace of all compact linear operators from $X$ into $Y$ by $\mathcal{K}(X,Y)$. If $X=Y$, then $\mathcal{B}(X,X) := \mathcal{B}(X)$ and $\mathcal{K}(X,X) := \mathcal{K}(X)$. For $T \in \mathcal{B}(X,Y)$, let $T'\in \mathcal{B}(Y^*,X^*)$ denote the Banach space adjoint of $T$ defined by $T'f=fT$ for all $f\in X^*$.
If $M$ is a closed subspace of Banach space $X$,
let $x+M=\{x+y:y\in M\}$ for $x\in X$,  
$X/M=\{x+M:x\in X\}$ and $\|x+M\|=\inf\{\|x+y\|:y\in M\}$. We know that $X/M$ with this norm is a Banach space. 
We let $M^\perp=\{f\in X^*:f(x)=0 \text{~for~}x\in M\}$. If $N$ is a closed subspace of $X^*$, then we let $^{\perp}\!N=\{x\in X:f(x)=0\text{~for~}f\in 
N\}$. Two closed subspaces $L,M$ of $X$ are said to be topologically complementary (or $L$ is topologically complemented by $M$) if the map $(x,y)\to x+y$ is an isomorphism from $L\times M$ (with the product topology) onto $X$.
Throughout the article, we let $T:X \to Y$ to be a nonzero bounded linear operator.

The spectral theorem for compact self-adjoint operators on a Hilbert space  states that if $(H,\left\langle \cdot,\cdot\right\rangle)$ is a Hilbert space and $T:H\to H$ is a compact and self-adjoint linear operator, then
\begin{align}\label{rep1}
	Tx=\sum_{n}\lambda_n\left\langle x,u_n\right\rangle u_n,\quad x\in H,
\end{align}
where $\lambda_1,\lambda_2,\ldots$ are the  eigenvalues of $T$ with $|\lambda_1|\geq|\lambda_2|\geq\cdots$. More generally, if $H_1,H_2$ are Hilbert spaces and $T:H_1\to H_2$ is a compact linear operator, then there exist orthogonal sets $\{u_n\}$ and $\{v_n\}$ of $H_1$ and $H_2$, respectively such that
	\begin{equation}\label{rep2}
	Tx=\sum_{n}\lambda_n\left\langle x,u_n\right\rangle v_n,\quad x\in H_1,
	\end{equation}
where $\lambda_1 \geq \lambda_2 \geq \cdots$
are the eigenvalues of the positive square root of $T^*T$, called the singular values of $T$
and  $u_1,u_2,\ldots$ are the corresponding eigenvectors which form an orthonormal basis for the orthogonal complement of $\ker(T)$. If $T$ is of finite rank, then $n$ in the above sums ranges from 1 upto $\rank(T)$. If $T$ is of infinite rank, then $\lambda_n\neq 0$ for all $n\in\mathbb{N}$ and $\displaystyle \lim_{n\to \infty}\lambda_n=0$. In this case, the series on the right side of the equation converges to $Tx$ and  if for $n\in\mathbb{N}$ and $x\in H$, $S_{n+1}(x)=\sum_{i=1}^n\left\langle x,u_i\right\rangle u_i$, then $\{TS_n\}$ converges to $T$ in the operator norm. 

The obvious difficulty with generalizing the above result for operators between Banach spaces is the absence of inner products in these spaces. However, if $H$ is a Hilbert space, then by the Riesz Representation Theorem, any $f:H\to \mathbb{F}$ is a bounded linear functional if and only if there exists $y \in H$ such that $f(x)= \langle x,y \rangle$ for all $x \in H$. By the Hahn-Banach Theorem, we know that for each nonzero element $x$ of a Banach space $X$, there exists $f\in X^*$ such that $\|f\|=1$ and $f(x)=\|x\|$.  If $X^*$  is strictly convex (its unit sphere does not contain any nontrivial line-segment), then such a functional is uniquely determined \cite{Taylor, Foguel}. Edmunds et al. used these facts to make the following generalization  of representation \eqref{rep2} for compact operators between reflexive Banach spaces with strictly convex dual spaces.
\begin{theorem}\cite{EEH2, EEH2_C}\label{Edmunds2}
	Let $X, Y$ be reflexive Banach spaces with strictly convex dual spaces and $T\in \mathcal{K}(X,Y)$. Then there exist sequences  $\{x_n\}$ in $S_X$ and $\{\xi_n\}$ in $X^*$ $($which are finite if $T$ is of finite rank and infinite otherwise$)$ with $\xi_i(x_j)=\delta_{ij}$ such that, if the sequence of operators $\{S_{n+1}\}$, given by $$S_{n+1}(x)=\sum_{i=1}^n\xi_i(x)x_i,\quad x\in X, n\in\mathbb{N}$$ is bounded in the operator norm, then
	\begin{align}\label{EEHrep1}Tx=\sum_{n}\xi_n(x)Tx_n,\quad x\in X,\end{align}	
	where $\displaystyle \lim_{n\to\infty}\|Tx_n\|=0$ if $T$ is of infinite rank.	
\end{theorem}

In \eqref{EEHrep1}, each $\xi_n$ is a certain linear combination of $f_1,\ldots, f_n$, where $f_n$ is the unique bounded linear functional  such that $\|f_n\|=1=\|x_n\|=f_n(x_n)$.
In Section 2 of the present paper, we prove the above results without assuming strict convexity of the dual spaces. To be precise, we prove the following. 
\begin{theorem}
	Let $X, Y$ be reflexive Banach spaces and $T\in \mathcal{K}(X,Y)$. Then there exist sequences  $\{x_n\}$ in $S_X$ and $\{\xi_n\}$ in $X^*$ $($which are finite if $T$ is of finite rank and infinite otherwise$)$ with $\xi_i(x_j)=\delta_{ij}$ such that, if the sequence of operators $\{S_{n+1}\}$, given by $$S_{n+1}(x)=\sum_{i=1}^n\xi_i(x)x_i,\quad x\in X, n\in\mathbb{N}$$ is bounded in the operator norm, then
	\begin{align}Tx=\sum_{n}\xi_n(x)Tx_n,\quad x\in X.\end{align}
	If $T$ is of infinite rank, then  $\displaystyle\lim_{n\to\infty}\|Tx_n\|=0$ and $\displaystyle\lim_{n\to\infty}\|TS_n-T\|=0$.
\end{theorem} 
  While in \cite{EEH2}, the fact that $\cap_nX_n= \ker(T)$, where  $X_{n+1}=\{x\in X:f(x)=0\text{~for~}f\in \lin\{f_1,\ldots,f_n\}\}$ for each $n\leq\rank(T)$, is proved as a consequence of Theorem \ref{Edmunds2}, we show the same without using the representation. In fact, this particular fact holds for any compact operator regardless of whether it satisfies the conditions for the existence of the representation.  In section 3, we discuss the representation for  a particular class of operators on Banach spaces, namely, those that \emph{behave like} compact self-adjoint or compact normal operators on Hilbert spaces. We show that for such operators, the representation has the same form as \eqref{rep1} and illustrate this with an example.

 We state a few fundamental results \cite{MF,Megginson} from  functional analysis, which we shall use in this paper. 
\begin{theorem}
	Let $T \in \mathcal{B}(X,Y)$ $(\text{resp.}~ \mathcal{K}(X,Y))$ and let $Z$ be a closed subspace of $X$. Define $\widetilde{T}:X/Z\to Y$ by 
	$\widetilde{T}(x+Z)=Tx$. Then $\widetilde{T} \in \mathcal{B}(X/Z,Y)$ $(\text{resp.}~ \mathcal{K}(X/Z,Y))$.
\end{theorem}

\begin{theorem}
	If $X, Y$ are reflexive Banach spaces and $T\in \mathcal{K}(X,Y)$, then $T$ attains its norm, that is, there exists $x \in S_X$ such that $\|Tx\| = \|T\|$.
\end{theorem}

\begin{theorem}\label{weakdense}
	Let $\{x_n\}$ be a sequence of a Banach space $X$ and $x\in X$. Then the following are equivalent.
	\begin{enumerate}
		\item $\{x_n\}$ converges weakly to $x$.
		\item $\{x_n\}$ is bounded and $\displaystyle\lim_{n\to\infty}f(x_n)=f(x)$ for all $f\in D$, where $D$ is a dense subset of $X^*$.
	\end{enumerate}
\end{theorem}

\begin{theorem}\label{perpMperp}
	If $M$ is a closed subspace of a Banach space $X$, then $(^{\perp}\!M)^\perp=M$.
\end{theorem}
\begin{theorem}\label{Dofquo}
	If $M$ is a closed subspace of a Banach space $X$, then $M^\perp$ is isometrically isomorphic to $(X/M)^*$ under the map $f\to\widetilde{f}$ with $\widetilde{f}(x+M)=f(x)$ for all $x\in X$.
\end{theorem}	
\begin{theorem}\label{complement}
	Two closed subspaces $L,M$ are topologically complementary if and only if $x\to x+M$ is an isomorphism from $L$ onto $X/M$.	
\end{theorem}

\section{Representation of compact operators}
Let $X$ be a Banach space. Given $x\in X$ and $x \neq 0$, let
$$\mathcal{J}^X_x=\{f\in S_{X^*}:f(x)=\|x\|\}.$$ The set $\mathcal{J}^X_x$ is nonempty by the Hahn-Banach Theorem. It can be easily verified that $\mathcal{J}^X_x$ is a closed, bounded and convex subset of $X$. The following result is pivotal in proving Edmunds et al.'s representation theorems for compact operators.

\begin{proposition}\label{DualitySM}\cite{EEH2}
Let $X, Y$ be reflexive Banach spaces with strictly convex dual spaces. Let $T \in \mathcal{K}(X,Y)$ and $a\in S_X$. The operator $T$ attains its norm at $a$ if and only if the unique members $f$ and $g$ of $\mathcal{J}^X_{a}$ and $\mathcal{J}^Y_{Ta}$, respectively satisfy the condition 
$$\frac{1}{\|T\|} T'g=\frac{1}{\|T\|}gT=f.$$
\end{proposition}

We generalize this result for arbitrary Banach spaces.
\begin{proposition}\label{Duality1}
	Let $X, Y$ be Banach spaces. Let $T \in \mathcal{B}(X,Y)$ and $a\in S_X$. Then the following are equivalent.
	\begin{enumerate}
		\item $\|Ta\|=\|T\|$. 
		\item $\frac{1}{\|T\|} T'g=\frac{1}{\|T\|}gT \in\mathcal{J}^X_a$ for all  $g\in\mathcal{J}^Y_{Ta}$.
		\item $\frac{1}{\|T\|} T'g=\frac{1}{\|T\|}gT\in\mathcal{J}^X_a$ for some  $g\in\mathcal{J}^Y_{Ta}$.
	\end{enumerate}
\end{proposition}
\begin{proof}
	(1)$\implies$ (2):	For any $g\in\mathcal{J}^Y_{Ta}$,
	$$\|T'g\|=\|gT\|\leq \|g\|\|T\|\leq \|T\|$$ and $$T'g(a)=gT(a)=g(Ta) =\|T(a)\|=\|T\|,$$ which proves that $\|T'g\|=\|T\|=\|T'\|$.
	Thus $\|\frac{1}{\|T\|}g\circ T\|=1$ and $\frac{1}{\|T\|}gT(a)=1=\|a\|$, which implies
	$\frac{1}{\|T\|} T'g=\frac{1}{\|T\|}gT\in\mathcal{J}^X_a$.
	
	(2)$\implies$(3): Trivial.
	
	(3)$\implies$(1): If $\frac{1}{\|T\|}gT\in\mathcal{J}^X_a$ for some  $g\in\mathcal{J}^Y_{Ta}$, then we have
	$$\frac{1}{\|T\|}g(Ta)=\|a\|\text{~and~} g(Ta)=\|Ta\|,$$
	which implies $\|Ta\|=\|T\|$.
\end{proof}

The following result is crucial in proving the representation theorem of compact operators.
\begin{lemma}\label{Nesting}
	Let $X, Y$ be  reflexive Banach spaces and $T \in \mathcal{K}(X,Y)$. Let $k\leq \rank(T)$. Then for each $j\in\{1,\ldots,k\}$, there exist
	\begin{enumerate}
		\item subspaces $X_j$ of $X$ and $Y_j$ of $Y$;
		\item $x_j\in S_{X_j}$; 
		\item $f_j\in\mathcal{J}^X_{x_j}$ and  $g_j\in\mathcal{J}^Y_{Tx_j}$
	\end{enumerate}
 such that
 \begin{enumerate}
 	\item [(a)]$X_1=X,Y_1=Y,  X_{j}=^{\perp}\!\!\lin\{f_1,\ldots,f_{j-1}\}, 
 	Y_{j}=^{\perp}\!\!\lin\{g_1,\ldots,g_{j-1}\}$ for $j>1$;
 	\item [(b)] $T_j(X_j)\subseteq Y_j$ and $\|Tx_j\| = \|T_j\| \neq 0$, where $T_j = T|_{X_j}$;
 	\item [(c)]$f_j|_{X_j}=\frac{1}{\|T_j\|}g_j T_j$ and $g_j|_{Y_j}\in \mathcal{J}^{Y_j}_{Tx_j}$.
 \end{enumerate}
\end{lemma} 
\begin{proof} We prove the result by induction on $k\in\mathbb{N}$ with $k\leq \rank(T)$.
	
	If $k=1$, then $X_1=X, Y_1=Y$ and $T_1=T\neq 0$. Since $T_1$ is norm attaining, there exists $x_1\in S_{X_1}$ such that $\|T_1\|=\|T_1 x_1\|=\|Tx_1\|$. Let $g_1\in\mathcal{J}^Y_{Tx_1}$ and  $f_1=\frac{1}{\|T_1\|}g_1 T_1\in \mathcal{J}^X_{x_1}$. 
	
	Suppose the statement holds for some $k\in\mathbb{N}$ and any $S\in \mathcal{K}(X,Y)$ such that $k\leq\rank(S)$. 
	Now, let $T \in \mathcal{K}(X,Y)$ be such that $k+1\leq \rank(T)$.  
	Let $X_1, Y_1,x_1, T_1, f_1,g_1$ 	 be as in the previous paragraph.  Let $X_2=^{\perp}\!\!\lin\{f_1\}, Y_2=^{\perp}\!\!\lin\{g_1\}$ and  $T_2=T|_{X_2}$. If $x\in X_2$, then by Proposition \ref{Duality1}, $\frac{1}{\|T_1\|}g_1 T_1(x)=f_1(x)=0$. 
	So we have $g_1(Tx)=0$, which proves that $Tx\in Y_2$. Thus $T(X_2) \subseteq Y_2$.
	
	Since for any $\alpha\neq 0$, we have $g_1(\alpha Tx_1)=\alpha\|T_1\|\neq 0$, it follows that $\lin\{Tx_1\}\cap Y_2=\{0\}$. Also, for any  $x\in X$, we have $f_1(x-f_1(x)x_1)=0$, which implies $x-f_1(x)x_1\in X_2$ and $Tx-f_1(x)Tx_1\in T(X_2)$. Hence $T(X)=\lin\{Tx_1\}\oplus T(X_2)$. It follows that
	$k\leq \rank(T_2)$. 
	By the induction assumption, for  $j\in\{1,\ldots,k\}$,  there exist
		\begin{enumerate}
		\item subspaces $X_{j+1}$ of $X_2$ and $Y_{j+1}$ of $Y_2$;
		\item $x_{j+1}\in S_{X_2}$;
		\item $\widehat{f}_{j+1}\in\mathcal{J}^{X_2}_{x_{j+1}}$ and $\widehat{g}_{j+1}\in\mathcal{J}^{Y_2}_{Tx_{j+1}}$
	\end{enumerate}
	such that
	\begin{enumerate}
		\item[(a)] $X_{j+1}=^{\perp}\!\!\lin\{\widehat{f}_2,\ldots,\widehat{f}_j\}$,
		$Y_{j+1}=^{\perp}\!\!\lin\{\widehat{g}_2,\ldots,\widehat{g}_j\}$ for $j>1$;
		\item [(b)] $T_{j+1}({X}_{j+1})\subseteq {Y}_{j+1}$ and $\|T_{j+1} x_{j+1}\|= \|T_{j+1}\| \neq 0$, where $T_{j+1} = T|_{X_{j+1}}$;
		\item[(c)] $\widehat{f}_{j+1}|_{X_{j+1}}=\frac{1}{\|T_{j+1}\|}\widehat{g}_{j+1} T_{j+1}$ and $\widehat{g}_{j+1}|_{Y_{j+1}}\in \mathcal{J}^{\widehat{Y}_{j+1}}_{Tx_{j+1}}$.
	\end{enumerate}
	
	For $j\in \{1,\ldots,k\}$, choose  $f_{j+1}\in\mathcal{J}^X_{x_{j+1}}$ and $g_{j+1}\in\mathcal{J}^Y_{Tx_{j+1}}$
	such that  $f_{j+1}|_{X_2}=\widehat{f}_{j+1}$ and $g_{j+1}|_{Y_2}=\widehat{g}_{j+1}$ by the Hahn-Banach Theorem. Now for $j\in \{2,\ldots,k\}$, we have  $X_{j+1}=X_2\cap^{\perp}\!\!\lin\{f_2,\ldots,f_j\}=^{\perp}\!\!\lin\{f_1\}\cap ^{\perp}\!\!\lin\{f_2,\ldots,f_j\}=^{\perp}\!\!\lin\{f_1\ldots,f_j\}$ and similarly, $Y_{j+1}=^{\perp}\!\!\lin\{g_1\ldots,g_j\}$. Thus for  $j\in\{1,\ldots,k+1\}$, we have
	\begin{enumerate}
		\item subspaces $X_j$ of $X$ and $Y_j$ of $Y$;
		\item $x_j\in S_{X_j}$; 
		\item $f_j\in\mathcal{J}^X_{x_j}$ and  $g_j\in\mathcal{J}^Y_{Tx_j}$
	\end{enumerate}
	such that
	\begin{enumerate}
		\item [(a)] $X_1=X,Y_1=Y,  X_{j}=^{\perp}\!\!\lin\{f_1,\ldots,f_{j-1}\}, 
		Y_{j}=^{\perp}\!\!\lin\{g_1,\ldots,g_{j-1}\}$ for $j>1$;
		\item [(b)]
		$T_j(X_j)\subseteq Y_j$ and $\|T_j x_j\|= \|T_j\| \neq 0,$ where $T_j =T|_{X_j}$;
		\item [(c)] $f_j|_{X_j}=\frac{1}{\|T_j\|}g_j T_j$ and $g_j|_{Y_j}\in \mathcal{J}^{Y_j}_{Tx_j}$.
	\end{enumerate}
This proves the statement for any $k\in\mathbb{N}$ and $T \in \mathcal{K}(X,Y)$ such that $k\leq\rank(T)$.
\end{proof}
\begin{theorem}\label{Nesting1}
Let $X, Y$ be  reflexive Banach spaces and $T \in \mathcal{K}(X,Y)$. Then for each $j\leq \rank(T)$, there exist
	\begin{enumerate}
		\item subspaces $X_j$ of $X$ and $Y_j$ of $Y$;
		\item $x_j\in S_{X_j}$;  
		\item $f_j\in\mathcal{J}^X_{x_j}$ and  $g_j\in\mathcal{J}^Y_{Tx_j}$
	\end{enumerate}
	such that
	\begin{enumerate}
		\item [(a)] $X_1=X,Y_1=Y,  X_{j}=^{\perp}\!\!\lin\{f_1,\ldots,f_{j-1}\}, 
		Y_{j}=^{\perp}\!\!\lin\{g_1,\ldots,g_{j-1}\}$ for $j>1$;
		\item [(b)]
		$T_j(X_j)\subseteq Y_j$ and $\|T_j x_j\|= \|T_j\| \neq 0,$ where $T_j =T|_{X_j}$;
		\item [(c)]$g_j|_{Y_j}\in \mathcal{J}^{Y_j}_{Tx_j}$ and $f_j|_{X_j}=\frac{1}{\|T_j\|}g_j T_j$.
	\end{enumerate}
\end{theorem}
\begin{proof}
	If $T$ if of rank $k$, then the statement holds by Lemma \ref{Nesting}. 
	
	If $T$ is of infinite rank, then we prove the theorem by the principle of recursive definition \cite{Munkres}. Let $X_1=X, Y_1=Y, T_1=T\neq 0$ and there exists $x_1\in S_{X_1}$ such that $\|T_1\|=\|T_1 x_1\|=\|Tx_1\|$. Let $g_1\in\mathcal{J}^Y_{Tx_1}$ and  $f_1=\frac{1}{\|T_1\|}g_1 T_1\in \mathcal{J}^X_{x_1}$. Given $n\in\mathbb{N}$, if for $j\in\{1,\ldots,n\}$ 
	\begin{enumerate}
		\item subspaces $X_j$ of $X$ and $Y_j$ of $Y$;
		\item $x_j\in S_{X_j}$; 
		\item $f_j\in\mathcal{J}^X_{x_j}$ and  $g_j\in\mathcal{J}^Y_{Tx_j}$
	\end{enumerate}
	are defined such that
	\begin{enumerate}
		\item [(a)] $X_1=X,Y_1=Y,  X_{j}=^{\perp}\!\!\lin\{f_1,\ldots,f_{j-1}\}, 
		Y_{j}=^{\perp}\!\!\lin\{g_1,\ldots,g_{j-1}\}$ for $j>1$;
		\item [(b)] $T_j(X_j)\subseteq Y_j$ and $\|T_j x_j\|= \|T_j\| \neq 0,$ where $T_j =T|_{X_j}$;
		\item [(c)] $f_j|_{X_j}=\frac{1}{\|T_j\|}g_j T_j$ and $g_j|_{Y_j}\in \mathcal{J}^{Y_j}_{Tx_j}$,
	\end{enumerate}
then let  $X_{n+1}=^{\perp}\!\!\lin\{f_1,\ldots,f_n\}$,
	$Y_{n+1}=^{\perp}\!\!\lin\{g_1,\ldots,g_n\}$
	and $T_{n+1}=T|_{X_{n+1}}$. There exists
	$x_{n+1}\in S_{X_{n+1}}$ such that $\|T_{n+1}\|=\|T_{n+1} x_{n+1}\|$. By the Hahn-Banach Theorem choose $f_{n+1}\in\mathcal{J}^X_{x_{n+1}}$ and  $g_{n+1}\in\mathcal{J}^Y_{Tx_{n+1}}$ such that $$f_{n+1}|_{X_{n+1}}=\frac{1}{\|T_{n+1}\|}g_{n+1} T_{n+1}\text{~and~} g_{n+1}|_{Y_{n+1}}\in \mathcal{J}^{Y_{n+1}}_{Tx_{n+1}}.\qedhere$$	
\end{proof}

Throughout the rest of this section, we let $X, Y$ be reflexive Banach spaces and $T \in \mathcal{K}(X,Y)$. For $j\leq \rank(T)$,  we let 
\begin{enumerate}
	\item subspaces $X_j$ of $X$ and $Y_j$ of $Y$ with $T_j=T|_{X_j}$;
	\item $x_j\in S_{X_j}$; 
	\item $f_j\in\mathcal{J}^X_{x_j}$ and  $g_j\in\mathcal{J}^Y_{Tx_j}$
\end{enumerate}
be as in Theorem \ref{Nesting1}. If $T$ is of rank $k$, let
$$X_{k+1}=^{\perp}\!\!\lin\{f_1,\ldots,f_k\}=\cap_{n=1}^{k+1} X_n.
$$ If $T$ is of infinite rank, let  $$X_{\infty}=^{\perp}\!\!\lin\{f_n:n\in\mathbb{N}\}=\cap_{n=1}^{\infty} X_n.$$

\begin{proposition}\label{linind1}
	If $1\leq j\leq \rank(T)$, then
	$\lin\{x_1,\ldots,x_j\}\cap X_{j+1}=\{0\}$ and $\lin\{Tx_1,\ldots,Tx_j\}\cap T(X_{j+1})=\{0\}$.
\end{proposition}
\begin{proof}
	Consider $\alpha_1,\ldots, \alpha_j\in\mathbb{F}$ such that
	$z=\alpha_1x_1+\cdots+\alpha_jx_j\neq0$. Let $i$ be the smallest in $\{1,\ldots,j\}$ such that $\alpha_i\neq 0$. Then 
	$f_i(z)=f_i(\alpha_ix_i)=\alpha_i\neq 0$ and $g_i(Tz)=g_i(\alpha_iTx_i)=\alpha_i\|T_i\|\neq0$, which imply $\lin\{x_1,\ldots,x_j\}\cap X_{j+1}=\{0\}$ and  $\lin\{Tx_1,\ldots,Tx_j\}\cap T(X_{j+1})=\{0\}$, respectively
\end{proof}

\begin{proposition}\label{Dofquo1}
	\begin{enumerate}
		\item If $T$ is of rank $k$, then $(X/X_{k+1})^*$ is isometrically isomorphic to $\lin\{f_1,\ldots, f_k\}$.
		
		\item If $T$ is of infinite rank, then $(X/X_{\infty})^*$ is isometrically isomorphic to $\overline{\lin}\{f_n:n\in \mathbb{N}\}.$
	\end{enumerate}
\end{proposition}
\begin{proof}
	By Theorems \ref{Dofquo} and \ref{perpMperp}, if $T$ is of rank $k$, then
	$(X/X_{k+1})^*$  is isometrically isomorphic to $X_{k+1}^{\perp}=\lin\{f_1,\ldots,f_n\}$.
  Similarly, if $T$ is of infinite rank, then $(X/X_{\infty})^*$ is isometrically isomorphic to $X_{\infty}^{\perp}$. By \cite[Lemma 15]{EEH1}, we have 
	$$X_{\infty}^{\perp}=\left(\cap_{n} X_{n+1}\right)^{\perp}= \overline{\cup_{n}X_{n+1}^{\perp}}=
	\overline{\cup_{n}\lin\{f_1,\ldots,f_n\}}=\overline{\lin}\{f_n:n\in\mathbb{N}\}.\qedhere
	$$ 
\end{proof}
\begin{proposition}\label{basis1}
	Suppose that $T$ is of rank $k$.  Let $x\in X$. If  $\alpha_1,\ldots,\alpha_k\in\mathbb{F}$ such that $x+X_{k+1}=\alpha_1x_1+\cdots+\alpha_kx_k+X_{k+1}$, then 
	$$\alpha_1=f_1(x),~~ \alpha_n=f_n(x)-\sum_{i=1}^{n-1}f_n(x_i)\alpha_i,\quad 1<n\leq k.$$
\end{proposition}
\begin{proof}
	We know that $f_i(x_i)=1$ and  $f_i(x_j)=0$ for $1\leq i<j\leq k$. 	By Proposition \ref{Dofquo1}, we have the following system of $k$ equations with $k$  unknowns $\alpha_1,\ldots,\alpha_k$.
	\begin{align*}
		\alpha_1&=f_1(x)\\
		\alpha_1f_2(x_1)+\alpha_2&=f_2(x)\\
		\vdots\\
		\alpha_1f_k(x_1)+\cdots+\alpha_{k-1}f_k(x_{k-1})+\alpha_k&=f_k(x)\\
	\end{align*}
	Solving the above system, we get 
	$$\alpha_1=f_1(x),~~ \alpha_n=f_n(x)-\sum_{i=1}^{n-1}f_n(x_i)\alpha_i,\quad 1<n\leq k.\qedhere$$
\end{proof}
\begin{proposition}\label{basis2}
	Suppose that $T$ is of infinite rank. Let $x\in X$.
	If $\{\alpha_n\}\subset \mathbb{F}$ such that $\sum_{i=1}^\infty\alpha_ix_i+ X_{\infty}$ converges weakly to $x+ X_{\infty}$, then $$\alpha_1=f_1(x),~~ \alpha_n=f_n(x)-\sum_{i=1}^{n-1}f_n(x_i)\alpha_i,\quad n>1.$$
\end{proposition}
\begin{proof}
	By Proposition \ref{Dofquo1},  we have  $$f(x)=\sum_{i=1}^\infty\alpha_if(x_i)$$
	for all $f\in\overline{\lin}\{f_n: n \in \mathbb{N}\}$. Since  $f_i(x_i)=1$ and  $f_i(x_j)=0$ for $i<j$, it follows that
	for a fixed  $n\in\mathbb{N}$,
	$$\sum_{i=1}^{j}\alpha_if_n(x_i)=\sum_{i=1}^{n}\alpha_if_n(x_i)\quad \text{~for~}\quad j\geq n.$$ So $$f_n(x)=\sum_{i=1}^\infty\alpha_if_n(x_i)=\sum_{i=1}^{n}\alpha_if_n(x_i).$$We prove by induction that the sequences in the statement are equal.
	Clearly, $\alpha_1=f_1(x)$. Assume that for some $n\in\mathbb{N}$ and for $k\in\{1\ldots,n\}$, $$f_k(x)-\sum_{i=1}^{k-1}f_k(x_i)\alpha_i=\alpha_k.$$ Then $$f_{n+1}(x)=\sum_{i=1}^{n+1}\alpha_if_{n+1}(x_i)=\sum_{i=1}^n\alpha_if_{n+1}(x_i)+\alpha_{n+1}.$$  Thus $$f_{n+1}(x)-\sum_{i=1}^n\alpha_if_{n+1}(x_i)=\alpha_{n+1}.$$  This completes the proof.
\end{proof}

Based on Propositions \ref{basis1} and \ref{basis2}, we define 
$$\xi_1=f_1, S_1=0$$
 and for $1\leq n\leq\rank(T)$,
$$\xi_{n+1}=f_{n+1}-\sum_{j=1}^{n}f_{n+1}(x_j)\xi_j,$$
$$
S_{n+1}x=\sum_{j=1}^{n}\xi_j(x)x_j,\quad x\in X.$$

\begin{proposition}\label{delta}
 For $1\leq i,j\leq \rank(T)$, $\xi_i(x_j)=\delta_{i,j}$.
\end{proposition}
\begin{proof}
We know that $f_i(x_i)=1$ and $f_i(x_j)=0$ for any $i,j\in \mathbb{N},i<j$.
So,  $\xi_1(x_j)=f_1(x_j)=\delta_{1,j}$. Suppose that for $i\in\{1,\ldots,n\}$ and $j\in\mathbb{N}$, $\xi_i(x_j)=\delta_{i,j}$.
Let $m\in\mathbb{N}$.\\
If $m>n+1$,
$$\xi_{n+1}(x_m)=f_{n+1}(x_m)-\sum_{i=1}^{n}f_{n+1}(x_i)\xi_i(x_m)=0-0=0.$$
If $m=n+1$, $$\xi_{n+1}(x_{n+1})=f_{n+1}(x_{n+1})-\sum_{i=1}^{n}f_{n+1}(x_i)\xi_i(x_{n+1})=1-0=1$$
If $m<n+1$,
$$\xi_{n+1}(x_m)=f_{n+1}(x_m)-\sum_{i=1}^{n}f_{n+1}(x_i)\xi_i(x_m)=f_{n+1}(x_m)-f_{n+1}(x_m)=0.\qedhere$$
\end{proof}
\begin{proposition}\label{funseq}
	Let $n\leq\rank(T)$. Then $X_{n+1} \subseteq \ker(S_{n+1})$ and $x-S_{n+1}(x)\in X_{n+1}$ for $x\in X$.
\end{proposition}
\begin{proof}
	For each $i\in\{1,\ldots,n\}$,  
	$$f_i(x-S_{n+1}(x))=f_i(x)-f_i\Big(\sum_{j=1}^{n}\xi_j(x)x_j\Big)=f_i(x)-\sum_{j=1}^{n}f_i(x_j)\xi_j(x).$$ For $i< j\leq k$, $f_i(x_j)=0$. So,
	\begin{align*}
		f_i(x-S_{n+1}(x))
		&=f_i(x)-\sum_{j=1}^{i}f_i(x_j)\xi_j(x)\\
		&=f_i(x)-\sum_{j=1}^{i-1}f_i(x_j)\xi_j(x)-\xi_i(x)\\&= \xi_i(x)-\xi_i(x)\\&=0.
	\end{align*}
	This proves that $x-S_{n+1}(x)\in X_{n+1}$.
\end{proof}
\begin{corollary}\label{directsum}
	For each $n\leq\rank(T)$, $X_{n+1} = \ker(S_{n+1})$ and
	$S_{n+1}$ is the linear projection of $X$ onto $\lin\{x_1,\ldots,x_n\}$.
\end{corollary}
\begin{proof} The proof follows from Proposition \ref{funseq}.\end{proof}

For each $n\in\mathbb{N}$ and $x\in X$, let $$\mathcal{P}_{n+1}(x)=\mathcal{P}_{X_{n+1}}(x)=\{z\in X_{n+1}:\|x-z\|=\|x+X_{n+1}\|\}.$$
\begin{theorem}\label{metpro}
	For each $n\in \mathbb{N}$,
	$$x-S_{n+1}(x)\in \mathcal{P}_{n+1}(x-S_n(x)).$$
\end{theorem}
\begin{proof}
	\begin{align*}
	\|(x-S_n(x))-(x-S_{n+1}(x))\|&=\|S_n(x)-S_{n+1}(x))\|\\
	&=\|\xi_n(x)x_n\|\\
	&=|\xi_n(x)|.
	\end{align*}
Let $y\in X_{n+1}$. Then  $\xi_i(y)=0$ for $i\in\{1,\ldots,n\}$.  
Also, by Proposition \ref{delta}, $\xi_n(S_n(x))=\sum_{i=1}^{n-1}\xi_i(x)\xi_n(x_i)=0$.  Since $x-S_n(x)\in X_n$, $\xi_i(x-S_n(x))=0$ for $i\in\{1,\ldots,n-1\}$. So we have,
\begin{align*}
|\xi_n(x)|&=|\xi_n(x-S_n(x)-y)|\\
&=|f_n(x-S_n(x)-y)-\sum_{i=1}^{n-1}f_n(x_i)\xi_i(x-S_n(x)-y)|\\
&=|f_n(x-S_n(x)-y)|\\
&\leq \|x-S_n(x)-y\|.
\end{align*}
This proves that $x-S_{n+1}(x)\in \mathcal{P}_{n+1}(x-S_n(x))$.
\end{proof}

\begin{theorem}\label{finranker}
	If $T$ is of rank $k$, then $X_{k+1}=\ker(T)$.
\end{theorem}
\begin{proof}
	By Proposition \ref{linind1}, $Tx_1,\ldots,Tx_k$ are linearly independent. So the range of $T$ must be equal to $\lin\{Tx_1,\ldots, Tx_k\}$. Also by Proposition \ref{linind1}, $\lin\{Tx_1,\ldots, Tx_k\}\cap T(X_{k+1})=\{0\}$. If $x\in X_{k+1}$, then $Tx\in \lin\{Tx_1,\ldots, Tx_k\}\cap T(X_{k+1})$. Hence $Tx=0$, which implies $X_{k+1}\subseteq \ker(T)$. 
	
	On the other hand, if $x\in \ker(T)$, then $Tx=0$. For each $i\in \{1,\ldots,k\}$, note that $g_iT_ix=g_iTx=0$ and  $f_i|_{X_i}=\frac{1}{\|T_i\|}g_i T_i$. Since $x\in X_1$, we have $$f_1(x)=\frac{1}{\|T_1\|}g_1 T_1(x)=0,$$ which in turn implies $x\in X_2$. Therefore, $f_2(x)=\frac{1}{\|T_2\|}g_2 T_2(x)$. Continuing in this manner, we can see that $f_1(x)=\cdots=f_k(x)=0$. Thus $x\in^{\perp}\!\!\lin\{f_1,\ldots,f_k\}=X_{k+1}$. This proves that $\ker(T)\subseteq X_{k+1}$.
\end{proof}

\begin{theorem}\label{FinRank}
	If $T$ is of rank $k$, then $\{x_n+\ker(T)\}_{n=1}^k$ is a basis for $X/\ker(T)$ with
	$$x+\ker(T)=\sum_{n=1}^k\xi_i(x)x_i+\ker(T),\quad x\in X$$ and
	$$Tx=\sum_{n=1}^k\xi_n(x)Tx_n,\quad x\in X.$$ 
	Further, if $\ker(T)$ is  topologically complemented by a subspace $Z$, then 
	$$x=\sum_{n=1}^k\xi_n(x)x_n,\quad x\in Z.$$
	
\end{theorem}
\begin{proof}
	By Theorem \ref{finranker}, $\ker(T)=X_{k+1}$. By Theorem \ref{funseq}, $x-S_{k+1}(x)\in X_{k+1}$. So for any $x\in X$, we have 
	\begin{align*}
		x+\ker(T)=x+X_{k+1}&=S_{k+1}(x)+X_{k+1}\\
		&=\sum_{i=1}^k\xi_i(x)x_i+X_{k+1}\\
		&=\sum_{i=1}^k\xi_i(x)x_i+\ker(T).
	\end{align*} 
	Since by Proposition \ref{finranker}, $T(X_{k+1})=0$, we have
	$T(x-S_{k+1}(x))=0$ for $x\in X$, which implies 
	$$Tx=\sum_{i=1}^k\xi_i(x)Tx_i,\quad x\in X.$$
If $\ker(T)$ is  topologically complemented by a subspace $Z$, then by Theorem \ref{complement},
$$x=\sum_{n=1}^k\xi_n(x)x_n,\quad x\in Z. \qedhere$$
\end{proof}

We now consider the case when $T$ is of infinite rank.

\begin{theorem}\label{limit}  
	If $T$ is  of infinite rank, then $\displaystyle \lim_{n\to\infty}\|T_n\|=0$.
\end{theorem}
\begin{proof}
	Let $T$ be of infinite rank. Then by Theorem \ref{Nesting1}, $T_n\neq 0$ for all $n$. For any $m<n$, we have $$\|T_m\|=\|Tx_m\|=g_m(Tx_m)=g_m(Tx_m-Tx_n)\leq \|Tx_m-Tx_n\|.$$ Since $T$ is a compact operator and $\{x_n\}$ is a bounded sequence, then $\{Tx_n\}$ has a convergent subsequence. It follows that $\{\|T_n\|\}$ has a subsequence converging to 0. Since $\{\|T_n\|\}$ is monotonically decreasing, we must have $\displaystyle\lim_{n\to\infty}\|T_n\|=0$.
	\end{proof}

\begin{theorem}\label{kernel}
	If $T$ is of infinite rank, then $\ker(T)=X_{\infty}$.
\end{theorem}
\begin{proof}
	 Let $x\in\ker(T)$. Since $X_{\infty}=^{\perp}\!\!\lin\{f_n:n\in\mathbb{N}\}=\cap_{n=1}^{\infty}X_n$, we use induction to prove that $x\in X_n$ for all $n\in\mathbb{N}$.  Clearly $x\in X_1$. Suppose that $x\in X_n$ for some $n\in\mathbb{N}$. Then  for each $i\in \{1,\ldots,n\}$, we have $x\in X_i$ and $f_i|_{X_i}=\frac{1}{\|T_i\|}g_iT$ from which it follows that  $f_i(x)=\frac{1}{\|T_i\|}g_i(Tx)=0$. Thus $x\in X_{n+1}=^{\perp}\!\!\lin\{f_1,\ldots,f_n\}$. This proves that $\ker(T)\subseteq X_{\infty}$.
	 
Since $X_{\infty}=^{\perp}\!\!\lin\{f_n:n\in\mathbb{N}\}=\cap_{n=1}^{\infty}X_n$, if $x\in\cap_{n=1}^{\infty}X_n$, then $\|Tx\|=\|T_n x\|\leq\|T_n\|\|x\| \to 0$ as $n \to \infty$, by Theorem \ref{limit}, which implies $\|Tx\|=0$. Hence $X_{\infty}\subseteq\ker(T)$.
\end{proof}

\begin{theorem}\label{weaklimit}
	For each $f\in\lin\{f_n:n\in \mathbb{N}\}$, $\displaystyle\lim_{n\to\infty}f(S_n(x))=f(x)$.
\end{theorem}
\begin{proof}
	Fix $k\in\mathbb{N}$, then for each $n\in\mathbb{N}$,
	$$f_k(S_{n+1}(x))
	=f_k\Big(\sum_{i=1}^n\xi_i(x)x_i\Big)
	=\sum_{i=1}^n\xi_i(x)f_k(x_i).$$
	If $n>k$, then $f_k(S_{n+1}(x))=\sum_{i=1}^k\xi_i(x)f_k(x_i)=f_k(S_{k+1}(x))$.
	Since $x-S_{k+1}(x)\in X_{k+1}=^{\perp}\!\!\lin\{f_1,\ldots,f_k\}$, we
	have $$f_k(x-S_{n+1}(x))=f_k(x-S_{k+1}(x))=0.$$ This proves that 
	$\displaystyle\lim_{n\to\infty}f_k(S_n(x))=f_k(x)$. It follows that  $\displaystyle \lim_{n\to\infty}f(S_n(x))=f(x)$
	for $f\in\cup_{k=1}^\infty\lin\{f_1,\ldots,f_k\}$.
\end{proof}

\begin{theorem}\label{boundedsn}
	 Let $T$ be of infinite rank and $x\in X$. If the sequence $\{S_n(x)+\ker(T)\}$ is bounded, then the following hold. 
	 \begin{enumerate}
	\item $\{S_n(x)+\ker(T)\}$ converges weakly to $x+\ker(T)$.
	\item $Tx=\sum_ {n=1}^\infty\xi_n(x)Tx_n$.
	\end{enumerate}
\end{theorem}
\begin{proof}
	By Theorem \ref{kernel}, $\ker(T)=X_{\infty}$.
 By Theorem \ref{weaklimit}, $\displaystyle\lim_{n\to\infty}f(S_n(x))=f(x)$ and so for $f\in\lin\{f_n:n\in\mathbb{N}\}$, $$\displaystyle \lim_{n\to\infty}\widetilde{f}(S_n(x)+\ker(T))=\lim_{n\to\infty}\widetilde{f}(S_n(x)+X_\infty)=\widetilde{f}(x+X_\infty)=\widetilde{f}(x+\ker(T)),$$ where $\widetilde{f}$ is the corresponding element in $(X/\ker(T))^*$ as in Theorem \ref{Dofquo} associated to $f$. Thus if $\{S_n(x)+\ker(T)\}$ is bounded, then by Theorems \ref{weakdense} and \ref{Dofquo1}, $\{S_n(x)+\ker(T)\}$ must converge weakly to $x+\ker(T)$. Since $\widetilde{T}$ is compact, $\{\widetilde{T}(S_n(x)+\ker(T))\}$  converges to $\widetilde{T}(x+\ker(T))$. But this is the same as saying $\{T(S_n(x))\}$ converges to $Tx$. 
	Thus $$Tx=\lim_{n\to\infty}\sum_{i=1}^n\xi_i(x)Tx_i=\sum_{n=1}^{\infty}\xi_n(x)Tx_n. \qedhere$$
\end{proof}
 \begin{theorem}\label{boundedsnu}
 	If $T$ is  of infinite rank and the sequence $\{S_n(x)+\ker(T)\}$ is bounded for all $x\in X$, then the following hold. 
 	 \begin{enumerate}
 	 	\item $\{x_n+\ker(T)\}$ is a basis for $X/\ker(T)$
 	 with $$x+\ker(T)=\sum_{n=1}^\infty\xi_n(x)x_n+\ker(T),\quad x\in X.$$
 	\item $Tx=\sum_ {n=1}^\infty\xi_n(x)Tx_n$ for  $x\in X$.
 	\item  $\{TS_n \}$ converges to $T$ is the operator norm.
 	\end{enumerate}
 	Further, if $\ker(T)$ is  topologically complemented by a subspace $Z$, then 
 $$x=\sum_{n=1}^\infty\xi_i(x)x_i,\quad x\in Z.$$
 \end{theorem}
\begin{proof}
	 (1) Let $x\in X$. By Theorem \ref{boundedsn}, the series
	$\sum_{n=1}^\infty \xi_n(x)x_n+\ker(T)$ converges weakly to $x+\ker(T)$.
	 By Proposition \ref{basis2}, if 
	 $\alpha_1,\alpha_2,\ldots\in\mathbb{F}$ such that $\sum_{n=1}^{\infty}\alpha_nx_n+\ker(T)$ converges weakly to $x+\ker(T)$, then $\alpha_n=\xi_n(x)$ for all $n$. 
	 Therefore, $\{x_n+\ker(T)\}$ is a weak basis for $X/\ker(T)$. By the Banach weak basis theorem \cite[Theorem 5.3]{Morrison}, $\{x_n+\ker(T)\}$ is a basis for $X/\ker(T)$.  
	 
	 (2) Proved in Theorem \ref{boundedsn}.
	 
	 (3) For $n\in\mathbb{N}$, define $\widehat{S}_n:X/\ker(T)\to X/\ker(T)$ by $$\widehat{S}_n(x+\ker(T))=S_n(x)+\ker(T).$$ For any $y\in\ker(T)$, we have
	 \begin{align*}
	\|\widehat{S}_n(x+\ker(T))\|&=\|S_n(x)+\ker(T)\|\\&=\|S_n(x+y)+\ker(T)\|\\&\leq
	\|S_n(x+y)\|\\&\leq\|S_n\|\|x+y\|,
	 \end{align*}
	 which gives
	  $\|\widehat{S}_n(x+\ker(T))\|\leq \|S_n\|\|x+\ker(T)\|$. Thus $\{\widehat{S}_n\}$ is a sequence of bounded linear operators. Since for each $x\in X$, $\{\widehat{S}_n(x+\ker(T))\}$ is a bounded sequence, it follows by the uniform boundedness principle that $\{\|\widehat{S}_n\|\}$ is bounded.  Also for $n\in\mathbb{N}$, $x-S_{n+1}(x)\in X_{n+1}$ for all $x\in X$, by Proposition \ref{funseq}. Now
	 \begin{align*}
	 	\|TS_{n+1}x-Tx\|&=\|T_{n+1}S_{n+1}x-T_{n+1}x\|\\
	 	&=\|\widetilde{T}_{n+1}(S_{n+1}x-x+\ker(T))\|\\
	 	&\leq \|\widetilde{T}_{n+1}\|[\|S_{n+1}x+\ker(T)\|+\|x+\ker(T)\|]\\
	 	&=\|T_{n+1}\|[\|\widehat{S}_{n+1}(x+\ker(T))\|+\|x+\ker(T)\|]\\
	 	&\leq \|T_{n+1}\|[\|\widehat{S}_{n+1}\|+1]\|x+\ker(T)\|\\
	 	&\leq \|T_{n+1}\|(\|\widehat{S}_{n+1}\|+1) \|x\|.
	 \end{align*}
 Thus $\|TS_{n+1}-T\|\leq \|T_{n+1}\|(\|\widehat{S}_{n+1}\|+1)$. Since $\{\|\widehat{S}_n\|\}$ is a bounded sequence and $\displaystyle\lim_{n\to\infty}\|T_{n+1}\|=0$, it follows that $\displaystyle\lim_{n\to\infty}\|TS_{n+1}-T\|=0$.
 \end{proof}

If we define $R_1=0$
and for $1\leq n\leq\rank(T)$, $$R_{n+1}f=\sum_{i=1}^nf(x_i)\xi_i,\quad f\in X^*,$$
then we have the following consequence of Theorems \ref{boundedsn} and \ref{boundedsnu}.
\begin{corollary}
	If $T$ is  of infinite rank and the sequence $\{S_n(x)+\ker(T)\}$ is bounded for all $x\in X$, then the following hold.
	\begin{enumerate}
		\item $T'g=\sum_nT'g(x_n)\xi_n,\quad g\in Y^*$.
		\item $\{R_{n+1}T'\}$ converges to $T'$ in the operator norm if $T$ is of infinite rank.
	\end{enumerate}
	
\end{corollary}
\begin{proof}
	(1) By Theorems \ref{boundedsn} and \ref{boundedsnu}, we know that 
	$$Tx=\sum_n\xi_n(x)Tx_n,\quad x\in X.$$
	For any $g\in Y^*$ and $x\in X$,
	$$T'g(x)=gTx=\sum_n\xi_n(x)gTx_n=\sum_ n\xi_n(x)T'g(x_n).$$
	Thus we have 
	$$T'g=\sum_nT'g(x_n)\xi_n,\quad g\in Y^*.$$
	
	(2) Let $n\in\mathbb{N}$, $g\in Y^*$ and $x\in X$.
	Then \begin{align*}
		\|R_{n+1}T'g(x)-T'g(x)\|&=\|R_{n+1}gT(x)-gT(x)\|\\
		&=\|\sum_{i=1}^ngT(x_i)\xi_i(x)-gT(x) \|\\
		&=\|g\left(\sum_{i=1}^nT(x_i)\xi_i(x)\right)-gT(x) \|\\
		&\leq \|g\|\left\|\sum_{i=1}^nT(x_i)\xi_i(x)-Tx \right\|\\
		&=\|g\|\|TS_{n+1}x-Tx\|\\
		&=\|g\|\|TS_{n+1}-T\|\|x\|.
		\end{align*}
	It follows that $\|R_{n+1}T'-T'\|\leq \|TS_{n+1}-T\|$.
	This proves that $\{R_{n+1}T'\}$ converges to $T'$ in the operator norm.
\end{proof}

\section{A special case}

Suppose that $H$ is a Hilbert space and $T \in \mathcal{B}(H)$ is positive, i.e., $T = T^*$ and $\langle Tx,x \rangle \geq 0$. It is known that $T$ is norm attaining if and only if $\|T\|$ is an eigenvalue of $T$ \cite[Corollary 2.4]{CN}. In particular, if $T$ is a positive compact operator, then for every invariant subspace $M$ of $H$, $T|_M$ attains its norm and therefore $\|T|_M\|$ is an eigenvalue of  $T|_M$. Recall that, a closed subpace $M$ of a Banach space $X$ is called an invariant subspace for $T \in \mathcal{B}(X,Y)$ if $T(M) \subseteq M$, i.e., $Tx \in M$ for all $x \in M$. If $T$ is a compact normal operator, then for every invariant subspace $M$ of $H$, $\gamma\|T|_M\|$ is an eigenvalue of  $T|_M$ for some $\gamma\in S_{\mathbb{F}}$.  In this section, we use the following classical fixed point theorem to obtain a representation for Banach space operators with this property.

\begin{theorem}$($Schauder Fixed Point Theorem$)$\label{SchauderFPT}\cite{Schauder, Zeidler} If $K$ is a closed bounded and convex subset of a Banach space $X$ and $f :K \to K$ is a continuous function that takes bounded sets to relatively compact, then $f$ has a fixed point in $K$.
\end{theorem}

\begin{theorem}\label{nonnegativeeigenvalue}
	Let  $T \in \mathcal{B}(X)$.
 If for every invariant subspace $M$ of $X$, $\|T|_M\|$ is an eigenvalue of $T|_M$, then every eigenvalue $\lambda$ of $T$ is nonnegative.
\end{theorem}
\begin{proof}
	Let $\lambda$ be an eigenvalue of $T$ and let $M$ be the null space of $T-\lambda I$. For each $x\in M$, $Tx=\lambda x\in M$. Thus $M$ is a nontrivial invariant subspace of $T$. Since
$\|T|_M\|=\sup\{\|Tx\|:x\in S_M\}=\sup\{|\lambda|\|x\|:x\in S_M\}=|\lambda|$ is an eigenvalue of $T|_M$, there exists $x_0\in M$, $x_0\neq 0$ such that $Tx_0=|\lambda|x_0$. But $Tx_0=\lambda x_0$ and so we must have $|\lambda|=\lambda$, which implies $\lambda\geq0$.
\end{proof}

\begin{corollary}\label{Injective}
	Let  $T \in \mathcal{B}(X)$ be such that for every 
	invariant subspace $M$ of $X$, $\|T|_M\|$ is an eigenvalue of $T|_M$. Then for any $\mu>0$,
	$\mu I+T$ is injective.
\end{corollary}
\begin{proof}
	Suppose there exists $x_0 \neq 0$ such that $(\mu I + T)x_0 = 0$. Then $Tx_0 =-\mu x_0$, which is a contradiction by Theorem \ref{nonnegativeeigenvalue}. \end{proof}

\begin{lemma}\label{scaling}
   Let $x\in X$ with $x \neq 0$. Then the following hold.
   \begin{enumerate}
   	\item $\mathcal{J}_{\gamma x}=\mathcal{J}_x$ for any $\gamma>0$.
   	\item $\mathcal{J}_{\gamma x}=\overline{\gamma}\mathcal{J}_{x}$ for any $\gamma \in S_{\mathbb{F}}$.
   \end{enumerate}
   \end{lemma}
\begin{proof}
\begin{enumerate}
	\item Let $\gamma>0$. If $f\in \mathcal{J}_{\gamma x}$, then $f(x)=\frac{1}{\gamma}f(\gamma x)=\frac{1}{\gamma}\gamma\|
	x\|=\|x\|$. So $\mathcal{J}_{\gamma x} \subseteq \mathcal{J}_x$.
	
	On the other hand, if $f\in\mathcal{J}_x$, then $f(\gamma x)=\gamma f(x)=\gamma\|x\|=\|\gamma x\|$, which implies $\mathcal{J}_{x} \subseteq \mathcal{J}_{\gamma x}$. Thus $\mathcal{J}_{\gamma x}=\mathcal{J}_x$.
	
	\item Let $\gamma\in S_{\mathbb{F}}$. If $f\in \mathcal{J}_{\gamma x}$, then $\|\gamma f\|=1$ and $\gamma f(x)=f(\gamma x)=\|\gamma x\|=\|x\|$, which proves that 
	$\gamma f\in \mathcal{J}_x$,  which implies $f\in \overline{\gamma}\mathcal{J}_x$. So $\mathcal{J}_{\gamma x} \subseteq \overline{\gamma}\mathcal{J}_x$.
	
	Next, if $f\in \overline{\gamma}\mathcal{J}_x$, then $\gamma f\in\mathcal{J}_x$. Therefore $\|f\|= \|\gamma f\|=1$ and $f(\gamma x)=\gamma f(x)=\|x\|=\|\gamma x\|$, which proves that $f\in \mathcal{J}_{\gamma x}$.
	Thus we have $\mathcal{J}_{\gamma x}=\overline{\gamma}\mathcal{J}_{x}.$
\end{enumerate}
\end{proof}

\begin{remark}\label{scalingforanyno}
	For any $\gamma \in \mathbb{F}\setminus \{0\}$, $\mathcal{J}_{\gamma x} = \frac{\overline{\gamma}}{|\gamma|} \mathcal{J}_x$.
\end{remark}

\begin{theorem}\label{Duality2}
	Let $X$ be a Banach space and $T \in \mathcal{B}(X)$. 
	Suppose there exists $a\in S_X$ and $\gamma\in S_{\mathbb{F}}$ such that $Ta=\gamma\|T\|a$. Then the following hold.
	\begin{enumerate}
		\item $\mathcal{J}_{Ta}=\overline{\gamma}\mathcal{J}_{a}$.
		\item 	$\frac{\overline{\gamma}}{\|T\|}T'$ maps $\mathcal{J}_{a}$ into itself.
	\end{enumerate} 
 
\end{theorem}
\begin{proof}
  By Remark \ref{scalingforanyno}, we have 
 $$\mathcal{J}_{Ta}=\mathcal{J}_{\gamma\|T\|a}=\overline{\gamma}\mathcal{J}_{a}.$$ Next, by Proposition \ref{Duality1}, we know that $\frac{1}{\|T\|}T'$ maps $\mathcal{J}_{Ta}$ into $\mathcal{J}_{a}$. So, $\frac{1}{\|T\|}T'$ maps $\overline{\gamma}\mathcal{J}_{a}$ into $\mathcal{J}_{a}$. If $f\in \mathcal{J}_{a}$,
then $$\frac{\overline{\gamma}}{\|T\|}T'(f)=\frac{1}{\|T\|}T'\left(\overline{\gamma}f\right)\in \mathcal{J}_{a}.\qedhere$$
\end{proof}

As a consequence, we have the following.
\begin{corollary}\label{FPT}
Let $X$ be a Banach space and $T \in \mathcal{K}(X)$.
If  $Ta=\gamma\|T\|a$ for some $a\in S_X$ and $\gamma\in S_{\mathbb{F}}$, then there exists $f\in\mathcal{J}_a$ such that $T'f=\gamma\|T\|f$.
\end{corollary}
\begin{proof}
	The proof follows by Theorems \ref{SchauderFPT} and \ref{Duality2}.
\end{proof}

	\begin{theorem}\label{Nesting2}
	Suppose $T \in \mathcal{K}(X)$ be such that for each
	invariant subspace $M$ of $X$, there exists $\gamma_M\in S_{\mathbb{F}}$ such that $\gamma_M\|T|_M\|$ is an eigenvalue of $T|_M$. Then for each $j\leq \rank(T)$, there exist
	\begin{enumerate}
		\item subspaces $X_j$ of $X$;
		\item $x_j\in S_{X_j}$,  $\gamma_j\in S_{\mathbb{F}}$, $f_j\in\mathcal{J}^X_{x_j}$
	\end{enumerate}
	such that
	\begin{enumerate}
		\item [(a)] $X_1=X,  X_{j}=^{\perp}\!\!\lin\{f_1,\ldots,f_{j-1}\}$ for $j>1$;
		\item [(b)] $T_j(X_j)\subseteq X_j$ and $Tx_j=\gamma_j\|T_j\|x_j$, where $T_j=T|_{X_j}$;
		\item [(c)] $f_j|_{X_j}=\frac{\overline{\gamma_j}}{\|T_j\|}f_jT_j$.
	\end{enumerate}
\end{theorem}

\begin{proof}
	The proof is similar to that of Theorem \ref{Nesting1} by using Corollary \ref{FPT}.
\end{proof}

In the next two results, we assume that $T\in \mathcal{K}(X)$ is such that for each
invariant subspace $M$ of $X$, there exists $\gamma_M\in S_{\mathbb{F}}$ such that $\gamma_M\|T|_M\|$ is an eigenvalue of $T|_M$. For each $j\leq \rank(T)$, let subspaces $X_j$ of $X$ with $T_j=T|_{X_j}$, $x_j\in S_{X_j}$, $\gamma_j\in S_{\mathbb{F}}$ and $f_j\in\mathcal{J}^X_{x_j}$ be as in Theorem \ref{Nesting2}.

\begin{theorem}
	 If $T$ is  of rank $k$, then 
$$Tx=\sum_{n=1}^k\lambda_n\xi_n(x)x_n,\quad x\in X,$$ where $\lambda_n=\gamma_n\|Tx_n\|=\gamma_n\|T_n\|$  with $|\gamma_n|=1$ for $n\in\{1,\ldots,k\}$ are eigenvalues of $T$ and $\xi_i(x_j)=\delta_{ij}$ for $i,j\in\{1,\ldots,k\}$. 
\end{theorem}
\begin{theorem}\label{boundedsnun}
		 If $T$ is of infinite rank and $\{S_n(x)+\ker(T)\}$ is bounded for all $x\in X$, then  
		$$Tx=\sum_ {n=1}^\infty\lambda_n\xi_n(x)x_n,\quad x\in X,$$ 
where $\lambda_n=\gamma_n\|Tx_n\|=\gamma_n\|T_n\|$  with $|\gamma_n|=1$ for $n\in\mathbb{N}$ are eigenvalues of $T$ and  $\lambda_n \to 0$ as $n \to \infty$.
\end{theorem}

\begin{example}
	For $x=(x_1,x_2,x_3,\ldots)\in\ell^2$, let
	$$\|x\|_{(k,1)}=\|(x_1,\ldots,x_k)\|_1+\|(x_{k+1},x_{k+2},\ldots)\|_2$$ and let $X=(\ell^2,\|\cdot\|_{(k,1)})$. Then $X^*=(\ell^2,\|\cdot\|_{(k,\infty)})$, where 	$$\|x\|_{(k,\infty)}=\max(\|(x_1,\ldots,x_k)\|_\infty,\|(x_{k+1},x_{k+2},\ldots)\|_2).$$
   Let $\{\alpha_n\}$ be a sequence of scalars converging to 0. Define $T:X\to X$ by $$T(x_1,x_2,x_3,\ldots)=(\alpha_1x_1,\alpha_2x_2,\alpha_3x_3,\ldots)\ \text{for}\ (x_1,x_2,x_3,\ldots) \in X.$$
\end{example}
The norm $\|\cdot\|_{(k,1)}$ on $X$ is equivalent to the standard norm $\|\cdot\|_2$. So the operator $T$ is a compact linear operator. Let $|\alpha_{n_1}|=\max\{|\alpha_n|:n\in\mathbb{N}\}$ and for $j\in\mathbb{N}$,
let $|\alpha_{n_{j+1}}|=\max\{|\alpha_n|:n\in\mathbb{N}\backslash\{n_1,\ldots,n_j\}\}$. 

Let $e_n(m) = \delta_{m,n}$ for all $m,n \in \mathbb{N}$.
Then for each $n\in\mathbb{N}$, we have $\mathcal{J}_{e_n}=\{x\in S_{X^*}:x_n=1\}$, which is not a singleton.

Let $X_1=X$ and $T_1=T$. Then $\|T_1\|=|\alpha_{n_1}|$ and $Te_{n_1}=|\alpha_{n_1}|e_{n_1}=\|T_1\|e_{n_1}$. Let $\alpha_{n_1}=\gamma_1\|T_1\|$. Define $f_1:=e_{n_1}\in\mathcal{J}_{e_{n_1}}$. Then for any $x\in X$, we have
$$T_1'f_1(x)=f_1T_1x=\left\langle T_1x,f_1\right\rangle=\left\langle T_1x,e_{n_1}\right\rangle=\alpha_{n_1}x_{n_1}=\gamma_1\|T_1\|f_1(x),$$ 
which proves that $f_1T_1=\gamma_1\|T_1\|f_1$.

 Let $X_2=^{\perp}\!\!\lin\{f_1\}=\{x\in X:\left\langle x,e_{n_1}\right\rangle=0\}=\{x\in X:x_{n_1}=0\}$ and $T_2=T|_{X_2}$. Then $T_2(X_2)\subset X_2$ and $\|T_2\|=|\alpha_{n_2}|$. Let $\alpha_{n_2}=\gamma_2\|T_2\|$. Define $f_2:=e_{n_2}\in \mathcal{J}_{e_{n_2}}$. Then 
 $$T_2'f_2(x)=f_2T_2x=\left\langle T_2x,f_2\right\rangle=\left\langle T_2x,e_{n_2}\right\rangle=\alpha_{n_2}x_{n_2}=\gamma_2\|T_2\|f_2(x) \ \text{for}\ x\in X_2.$$ Thus $f_2T_2=\gamma_2\|T_2\|f_2|_{X_2}$.

  Proceeding in this manner, we obtain the following for $j\leq \rank(T)$:
 \begin{enumerate}
 \item[] $X_j=\{x\in X:x_{n_1}=\cdots=x_{n_{j-1}}=0\}$; 
 \item[] $T_j=T|_{X_j}$;
 \item[] $Te_{n_j}=\alpha_{n_j}e_{n_j}=\gamma_j\|T_j\|e_{n_j}$;
 \item[] $f_j=e_{n_j}\in\mathcal{J}_{e_{n_j}}$
 such that $f_j|_{X_j}=\frac{\overline{\gamma_j}}{\|T_j\|}f_jT_j$.
 \end{enumerate}

Further, $\xi_1=f_1$ and for $j\in\mathbb{N}$,
 and $$\xi_{j+1}=f_{j+1}-\sum_{i=1}^j f_{j+1}(e_{n_i})\xi_i=f_{j+1}-\sum_{i=1}^j \left\langle e_{n_i},e_{n_{j+1}}\right\rangle\xi_i=f_{j+1}.$$
So $j\in\mathbb{N}$ and for $x\in X$, we get
 $$S_{j+1}(x)=\sum_{i=1}^j\xi_i(x)e_{n_i}=\sum_{i=1}^jf_i(x)e_{n_i}=
 \sum_{i=1}^j\left\langle x,e_{n_i}\right\rangle e_{n_i}=\sum_{i=1}^jx_{n_i} e_{n_i}.$$
 Clearly, $\|S_{j+1}(x)\|\leq \|x\|_{(k,1)}$ for all $x\in X$, which proves that $\|S_{j+1}\|\leq 1$ for $j\in \mathbb{N}$.
 
 Hence by Theorems \ref{FinRank} or \ref{boundedsnun} depending on whether $T$ is of finite or infinite rank, we have 
 $$Tx=\sum_ j\alpha_{n_j}x_{n_j}e_{n_j},\quad x\in X$$
 and $$T'y=\sum_{j}\alpha_{n_j} y_{n_j}e_{n_j},\quad y\in X^*.$$
 

\section*{Acknowledgement}
The research of the first two authors is supported by MTR/2019/001307, SERB, Department of
Science and Technology (DST), Government of India. The research of the third author is supported by the DST-INSPIRE Fellowship (Grant No. DST/INSPIRE FELLOWSHIP/2018/IF180107), Government of India.

\bibliographystyle{amsplain}

\end{document}